\documentclass[11pt]{article}
\usepackage{amsmath,amssymb}
\usepackage{mathrsfs}
\usepackage{tikz}
\usetikzlibrary{calc}

\newtheorem{thm}{Theorem}[section]
\newtheorem{lem}[thm]{Lemma}

\newcommand{\p}{\mathcal{P}}
\newcommand{\e}{\varepsilon}
\newcommand{\T}{\mathcal{T}}
\newcommand{\C}{\mathrm{cr}}

\def\pf{\noindent{\it Proof.} }
\def\qed{\nopagebreak\hfill{\rule{4pt}{7pt}}\medbreak}

\makeatletter \@addtoreset{equation}{section} \makeatother

\makeatletter
\newcommand{\rn}[1]{($\romannumeral #1$)}
\newcommand{\rnp}[1]{($\romannumeral #1'$)}
\makeatother

\begin{document}

\begin{center}
{\Large\bf Determining All Universal Tilers}
\end{center}

\begin{center}
David G. L. Wang\\[6pt]

Beijing International Center for Mathematical Research\\
Peking University, Beijing 100871, P. R. China\\
{\tt wgl@math.pku.edu.cn}
\end{center}

\begin{abstract}
A universal tiler is a convex polyhedron
whose every cross-section tiles the plane.
In this paper, we introduce a certain slight-rotating operation
for cross-sections of pentahedra.
Based on a selected initial cross-section
and by applying the slight-rotating operation suitably,
we prove that a convex polyhedron is a universal tiler
if and only if it is a tetrahedron
or a triangular prism.
\end{abstract}


\noindent\textbf{2010 AMS Classification:} 05B45 52C20

%05B45 Tessellation and tiling problems
%52C20 Tilings in 2 dimensions

\section{Introduction}

A tiler is a polygon that can cover the plane
by congruent repetitions without gaps or overlaps.
The problem of determining all tilers,
called alternatively the problem of tessellation
or plane tiling, is one of the most famous problems
in discrete mathematics,
and is still open to the best of our knowledge.
For a whole theory of tessellation,
see Grunbaum and Shephard's book~\cite{GS87}.

Considering a variant of the problem of plane tiling,
Akiyama \cite{Aki07} found all convex polyhedra
whose every development is a tiler.
The key idea in his proof is to find
whether the facets of a polyhedron
tile the plane in certain stamping manner.
Noticing that facets are special cross-sections,
we study another variant of the problem of plane tiling
that what kind of polyhedra are so well-performed that
every its cross-section is a tiler.

Let~$\p$ be a convex polyhedron,
and~$\pi$ a plane.
Denote the intersection of~$\pi$ and~$\p$ by~$C(\pi)$.
We say that~$\pi$ intersects~$\p$ non-trivially if
$C(\pi)$ is a non-degenerated polygon,
that is, $C(\pi)$ has at least~$3$ edges.
We call~$C(\pi)$ a {\em cross-section}
if~$\pi$ crosses~$\p$ nontrivially.
The polyhedron~$\p$ is said to be a {\em universal tiler}
if every cross-section of~$\p$ is a tiler.
In this paper, we will determine all universal tilers.

A triangular prism is a pentahedron with parallel facets.
With the aid of Euler's formula,
and Reinhardt's theorem~\cite{Rei18}
for the results on tilers with $n$ $(n\ne5)$ edges,
the author~\cite{Wang11X} managed to obtain
the following necessary condition
for the number of faces of a universal tiler
by suitably choosing cross-sections
of a given polyhedron.

\begin{thm}
A convex polyhedron is a universal tiler
only if it is a tetrahedron or a pentahedron.
Moreover, every tetrahedron and every triangular prism
is a universal tiler.
\end{thm}

In light of the above theorem,
the problem of determining all universal tilers
turns out to be the one of finding the list of
pentahedron universal tilers.
One of the difficulties in determining whether a pentahedron
is a universal tiler is the fact that no one knows
the list of pentagonal tilers,
although there are $14$ classes of pentagonal tilers
are found, see Hirschhorn and Hunt~\cite{HH85},
and Sugimoto and Ogawa~\cite{SO06} for instance.

The key idea used in solving the universal tiler problem
consists of two parts.
One is to select an initial cross-section from a given pentahedron
subject to some technical conditions.
It is an extension of the method adopted in~\cite{Wang11X}.
The other is to suitably apply a certain {\em slight-rotating operation}
based on the initial cross-section.
By suitably applying the operation at most three times,
and considering the local situations of the tessellations,
we can prove that only triangular prisms are pentahedron universal tilers.
The whole proof has nothing to do with
the knowledge of the complete list of pentagonal tilers.

This paper is organized as follows.
In the next section, we give necessary notion and notation
on tessellations of the plane by a single polygon.
Section 3 is devoted to select the initial cross-section subject to
some technical conditions and to introduce the slight-rotating operation.
In the last section,
by applying the slight-rotating operation
we prove that one may invariantly obtain a non-tiler
cross-section of a pentahedron~$\p$ unless~$\p$
is a triangular prism.

\section{Preliminary}

In this section, we introduce some necessary notion and notation.
Suppose that
\[
T=V^1V^2\cdots V^5
\]
is a pentagonal tiler.
Let~$\mathcal{T}$ be a tessellation of the plane by copies of~$T$.
Denote the copies used in~$\mathcal{T}$ by
\[
\{\,T_i=V_i^1V_i^2\cdots V_i^5\,\colon\ i\in\Lambda\,\},
\]
where~$\Lambda$ is a set.
Then every~$T_i$ has the same shape as~$T$.

Let~$i\in\Lambda$ and~$\e>0$.
Since~$T$ is a tiler,
the $\e$-neighborhood of the point~$V_i^j$ in~$\T$
must be covered without gaps or overlaps.
It follows that either there is a sequence
\[
V_i^j,\ \ V_{i_1}^{j_1},\ \ V_{i_2}^{j_2},\ \ \ldots,\ \ V_{i_k}^{j_k}
\]
($k\ge2$) of angles arranged counter-clockwise
which fulfilled the whole $2\pi$-area around the point $V_i^j$
(see the left part of Figure 1),
or there is a sequence
\[
V_{i_1}^{j_1},\ \ \ldots,\ \ V_{i_k}^{j_k},\ \
V_i^j,\ \ V_{s_1}^{t_1},\ \ \ldots,\ \ V_{s_h}^{t_h},
\]
($k+h\ge1$) of angles arranged counter-clockwise which fulfilled
a $\pi$-area around $V_i^j$
(see the right part of Figure 1).
In the former case, we denote the local tessellation
around the point~$V_i^j$ by
\begin{equation}\label{S1}
S(V_i^j)
=\bigl[\,V_i^j,\ V_{i_1}^{j_1},\ V_{i_2}^{j_2},\ \ldots,\ V_{i_k}^{j_k}\,\bigr].
\end{equation}
In the latter case, we denote
\begin{equation}\label{S2}
S(V_i^j)
=\bigl[\,V_i^j,\ V_{s_1}^{t_1},\ \ldots,\ V_{s_h}^{t_h},\ \pi,\
V_{i_1}^{j_1},\ \ldots,\ V_{i_k}^{j_k}\,\bigr].
\end{equation}
\begin{center}
\begin{tikzpicture}
\begin{scope}
\draw
(0,0)node[above=4mm,right=2mm]{$V_i^j$}--(1.5,0)
(0,0)node[above=3mm]{$V_{i_1}^{j_1}$}--(65:1.3)
(0,0)node[left=2mm]{$V_{i_2}^{j_2}$}--(130:1.5)
(0,0)node[below=4mm]{$\cdots$}--(220:1.3)
(0,0)node[below=4mm,right=2mm]{$V_{i_k}^{j_k}$}--(300:1);

\fill[black,opacity=1] (0,0) circle (1.5pt);
\end{scope}

\begin{scope}[xshift=6cm]
\draw
(0,-.6)node[above=3.7mm,right=15mm]{$V_{i_1}^{j_1}$}--(4,-.6)
(0,-.6)node[above=8mm,right=10mm]{$\cdots$}--+(20:3.5)
(0,-.6)node[above=12mm,right=3mm]{$V_{i_k}^{j_k}$}--+(40:3)
(0,-.6)node[above=10mm,left=-3mm]{$V_i^j$}--+(75:2)
(0,-.6)node[above=14mm,left=6mm]{$V_{s_1}^{t_1}$}--+(115:2)
(0,-.6)node[above=9mm,left=12mm]{$\cdots$}--+(140:2.5)
(0,-.6)node[above=3.7mm,left=15mm]{$V_{s_h}^{t_h}$}--+(160:3)
(0,-.6)node[above=4mm,left=2mm]{}--(-3,-.6);

\fill[black,opacity=1] (0,-.6) circle (1.5pt);
\end{scope}

\begin{scope}[xshift=28mm,yshift=-16mm]
\node[xshift=15mm,text centered]
{Figure 1. The local tessellation around the point $V_i^j$.};
\end{scope}
\end{tikzpicture}
\end{center}

Suppose that there is a sequence
\[
V_{i_1}^{j_1}V_{i_1}^{k_1},\quad
V_{i_2}^{j_2}V_{i_2}^{k_2},\quad
\ldots,\quad
V_{i_s}^{j_s}V_{i_s}^{k_s}
\]
of $s$ ($s\ge1$) edges
such that the copies
\[
T_{i_1},\ T_{i_2},\ \ldots,\ T_{i_s}
\]
lie on the same side of the line segment
\[
\vec{l}=V_{i_1}^{j_1}V_{i_s}^{k_s}.
\]
To be more precise,
the point~$V_{i_r}^{k_r}$
coincides with the point~$V_{i_{r+1}}^{j_{r+1}}$
for each $1\le r\le s-1$.
In this case, we write
\[
\vec{l}=
V_{i_1}^{j_1}V_{i_1}^{k_1}+
V_{i_2}^{j_2}V_{i_2}^{k_2}+\cdots+
V_{i_s}^{j_s}V_{i_s}^{k_s}
\]
Assume that there also holds
\[
\vec{l}=
V_{i_1'}^{j_1'}V_{i_1'}^{k_1'}+
V_{i_2'}^{j_2'}V_{i_2'}^{k_2'}+\cdots+
V_{i_t'}^{j_t'}V_{i_t'}^{k_t'}
\]
such that all the copies
\[
T_{i_1'},\ T_{i_2'},\ \ldots,\ T_{i_t'}
\]
lie on the other side of $\vec{l}$,
as illustrated in Figure 2.
\begin{center}
\begin{tikzpicture}
\begin{scope}
\draw
(0,0)coordinate(V1)--(8,0)
(0,0)node[above]{$V_{i_1}^{j_1}$}--(150:1.2)
(2,0)node[above=4mm, left=0mm]{$V_{i_1}^{k_1}$}
coordinate(i1)--(1.8,1)
(2,0)node[above=4mm, right=2mm]{$V_{i_2}^{j_2}$}--(2.5,1)
(4,0)node[above=4mm, left=-2mm]{$V_{i_2}^{k_2}$}
coordinate(i2)--(4.2,1)
(i2)node[above=4mm, right=6mm]{$\cdots$}
(6,0)node[above=4mm, right=-2mm]{$V_{i_s}^{j_s}$}
coordinate(js)--(5.7,1)
(8,0)node[above=4mm, left=-3mm]{$V_{i_s}^{k_s}$}
coordinate(ks)--(9,1)

(0,0)node[below]{$V_{i_1'}^{j_1'}$}--(200:1.2)
(2.5,0)node[below=5mm, left=-1mm]{$V_{i_1'}^{k_1'}$}
coordinate(i1')--(2.6,-1)
(i1')node[below=5mm, right=8mm]{$\cdots$}
(5,0)node[below=5mm, right=-2mm]{$V_{i_t'}^{j_t'}$}
coordinate(jt')--(4.5,-1)
(8,0)node[below=5mm, left=-2mm]{$V_{i_t'}^{k_t'}$}
coordinate(kt')--(9,-1);

\foreach \point in {V1,i1,i2,js,ks,i1',jt',kt'}
\fill[black,opacity=1] (\point) circle (1.5pt);
\end{scope}

\begin{scope}[xshift=25mm,yshift=-16mm]
\node[xshift=15mm,text centered]
{Figure 2. The arrangement~\eqref{S-e}.};
\end{scope}
\end{tikzpicture}
\end{center}
In this case, we say that~$\vec{l}$
is represented in~$\mathcal{T}$, denoted as
\begin{equation}\label{S-e}
S\Bigl(\,V_{i_1}^{j_1}V_{i_1}^{k_1},\ \
V_{i_2}^{j_2}V_{i_2}^{k_2},\ \
\ldots,\ \
V_{i_s}^{j_s}V_{i_s}^{k_s}\,\Bigr)
=\Bigl[\,V_{i_1'}^{j_1'}V_{i_1'}^{k_1'},\ \
V_{i_2'}^{j_2'}V_{i_2'}^{k_2'},\ \
\ldots,\ \
V_{i_t'}^{j_t'}V_{i_t'}^{k_t'}\,\Bigr].
\end{equation}

Denote by~$\mathbb{Z}$ the set of integers,
by~$\mathbb{Z}^+$ the set of positive integers,
and by~$\mathbb{N}$ the set of nonnegative integers.
Let~$N>0$ and~$a_1,\ldots,a_5\in\mathbb{N}$.
We say that the set
\begin{equation}\label{sum-N-collection}
\{\,a_1V^1,\
a_2V^2,\
\ldots,\
a_5V^5\,\}
\end{equation}
is a sum-$N$-collection of the angles if
\[
N=\sum_{i=1}^5a_iV^i.
\]
We call~\eqref{sum-N-collection} a sum-$N$-collection to the angle~$V^j$
if~$a_j\ge1$.
For convenience,
we remove the term~$a_iV^i$
from the collection~\eqref{sum-N-collection}
if
\[
a_i=0.
\]
Denote by~$R_N(V^j;\,T)$ the set of sum-$N$-collections to~$V^j$ in~$T$.
For simplifying notation, we write
\[
R_{2\pi}(V^j;\,T)=R\,(V^j;\,T).
\]
Since $T$ is a tiler, we have
\[
R\,(V^j;\,T)\ne\varnothing
\]
for any~$1\le j\le 5$.
For example, if
\begin{align*}
&V^1=V^2={5\pi\over 6},\\[5pt]
&V^3={2\pi\over 3},\\[5pt]
&V^4=V^5={\pi\over 3},
\end{align*}
then we have
\begin{align*}
R\,(V^1;\,T)&=\bigl\{
\{2V^1,\,V^4\},\
\{2V^1,\,V^5\},\
\{V^1,\,V^2,\,V^4\},\
\{V^1,\,V^2,\,V^5\}\bigr\},\\[5pt]
R_\pi(V^3;\,T)&=\bigl\{
\{V^3,\,V^4\},\
\{V^3,\,V^5\}\bigr\}.
\end{align*}

Note that the sum of any four angles of a pentagon
is larger than~$2\pi$, while the sum of any three angles
is larger than~$\pi$. This leads to the following lemma immediately.

\begin{lem}\label{lem-sum-collection}
Any sum-$2\pi$-collection~\eqref{sum-N-collection}
has at most three positive~$a_j$,
while any sum-$\pi$-collection~\eqref{sum-N-collection}
has at most two positive~$a_j$.
\end{lem}

\section{The slight-rotating operation}

In this section,
we first demonstrate the selection of the initial cross-section
subject to some technical conditions.
With the initial cross-section in hand,
we can then introduce the slight-rotating operation
which will play the essential role in determining all pentahedron universal tilers.

Denote by~$\mathscr{E}$ the set of pentahedron universal tilers
without parallel facets.
The goal of this paper is to prove that
\[
\mathscr{E}=\varnothing.
\]
We will do this by contradiction.
It is well-known that pentahedra have two distinct topological types.
One is the quadrilateral-based pyramids,
the other is composed of
a pair of triangular bases
and three quadrilateral sides,
such as triangular prisms.
In particular, we see that any pentahedron has a quadrilateral facet.
Throughout this paper, we suppose that
$\p\in\mathscr{E}$,
and
\[
Q=ABCD
\]
is a convex quadrilateral facet of~$\p$.

Let~$E$ be a point lying in the interior of the line segment~$AB$
such that $E$ is neither~$A$ nor~$B$.
In this case, we express
\[
E\in AB.
\]
Set another point~$F\in BC$.
\begin{center}
\begin{tikzpicture}
\draw(0,0)node[below]{$A$}coordinate(A)
--(3,0)node[below]{$B$}coordinate(B)
--(4,1)node[right]{$C$}coordinate(C)
--(1,1.5)node[left]{$D$}
--cycle;
\coordinate (E) at ($(A)!0.2!(B)$);
\coordinate (F) at ($(B)!0.2!(C)$);
\draw[thin] (E)node[below]{$E$}--(F)node[right]{$F$};
\node[above of=E,yshift=-8mm,xshift=2mm]{};
\node[above of=F,yshift=-8mm,xshift=-1mm]{};
\foreach \point in {E,F}
\fill[black,opacity=1](\point) circle (1.5pt);
\begin{scope}[yshift=-12mm]
\node[xshift=2cm,text centered]
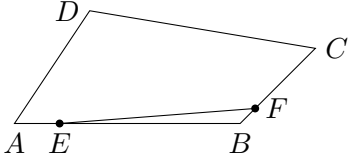
{Figure 3. The pentagon~$AEFCD$.};
\end{scope}
\end{tikzpicture}
\end{center}
As illustrated in Figure 3, we have a convex pentagon
\[
T=AEFCD.
\]
The method of finding a cross-section
by choosing $E\to A$ and $F\to B$
was used in~\cite{Wang11X}, however,
at this time, we need to select the points~$E$ and~$F$ more carefully.
Before formulating the conditions for selecting~$E$ and~$F$,
we give the following lemma which will be frequently used
in the proof of Lemma~\ref{lem-EF}.

\begin{lem}\label{lem-finite}
Let~$N\ge0$. Then both the set
\[
\bigl\{\theta>B/2\,\colon\
N=a\!\cdot\!\theta+b\!\cdot\!\phi+c\!\cdot\!\psi,\
a\in\mathbb{Z}^+,\
b,c\in\mathbb{N},\
\phi,\psi\in\{A,C,D\}\bigr\}
\]
and the set
\[
\bigl\{l<N\,\colon\
l=a\!\cdot\!CD+b\!\cdot\!DA+c\!\cdot\!AE,\
a,b\in\mathbb{N},\,c\in\{0,\pm1,\pm2\}\bigr\}
\]
have finite cardinalities.
\end{lem}

It is easy to prove the above lemma and we omit the proof.

\begin{lem}\label{lem-EF}
There exists~$\delta_E>0$ such that
for any point~$E\in AB$ with~$AE<\delta_E$,
there exists~$\delta_F=\delta_F(E)>0$ such that
for any point~$F\in BC$ with~$BF<\delta_F$,
the pentagon
\[
T=AEFCD
\]
satisfies
\begin{itemize}
\item[\rn{1}]
$2AE<\min\{EF,\,FC,\,CD,\,DA\}$;
\item[\rn{2}]
$R\,(E;\,T)=R\,(F;\,T)\subseteq
\{\{E,\,F,\,c\gamma\}\,\colon\
c\in\mathbb{Z}^+,\
\gamma\in\{A,C,D\}\}$.
\item[\rn{3}]
for any~$a,b\in\mathbb{N}$, $c\in\{0,\,\pm1,\,\pm2\}$, and~$l\in\{EF,\,2FC\}$,
there holds
\[
l\ne a\!\cdot\!CD+b\!\cdot\!DA+c\!\cdot\!AE.
\]
\end{itemize}
\end{lem}

\pf Choose a point~$E\in AB$ and move it towards~$A$
such that the point~$E$ can be arbitrarily closed to
but never arrive at the point~$A$.
Similarly, choose~$F\in BC$ and
move it towards~$B$
such that~$F$ can be arbitrarily closed to
but never arrive at the point~$B$.
Let
\[
x=\max\{AE,\,BF\}.
\]
In the above moving procedure, it is clear that
\begin{align}
&\lim_{x\to0}AE=0,\label{lim-AE}\\[5pt]
&\lim_{x\to0}EF=AB,\label{lim-EF}\\[5pt]
&\lim_{x\to0}FC=BC,\label{lim-FC}\\[5pt]
&\lim_{x\to0}\angle AEF=\pi,\label{lim-angle-E}\\[5pt]
&\lim_{x\to0}\angle EFC=B.\label{lim-angle-F}
\end{align}
By~\eqref{lim-AE}, \eqref{lim-EF} and~\eqref{lim-FC}, there exists
$\delta_1$ such that the condition~\rn{1} holds for any~$x<\delta_1$.

By Lemma~\ref{lem-sum-collection}, we have
\begin{equation}\label{eq2}
R\,(E;\,T)\subseteq
\bigl\{\{aE,\,b\beta,\,c\gamma\}\,\colon\
a\in\mathbb{Z}^+,\ b,c\in\mathbb{N},\ \beta,\gamma\in\{F,A,C,D\}\bigr\}.
\end{equation}
By~\eqref{lim-angle-E} and~\eqref{lim-angle-F},
there exists~$\delta_2<\delta_1$
such that for any~$x<\delta_2$, there holds
\[
2E+\alpha>2\pi,\quad \forall\,\alpha\in\{A,F,C,D\}.
\]
Let~$x<\delta_2$. Then~$a=1$ in~\eqref{eq2}.
By Lemma~\ref{lem-finite} and the limit~\eqref{lim-angle-E},
there exists~$\delta_3<\delta_2$ such that
for any~$x<\delta_3$, we have
\[
E+b\beta+c\gamma\ne2\pi,\quad
\forall\, b,c\in\mathbb{N},\
\beta,\gamma\in\{A,C,D\}.
\]
Let~$x<\delta_3$.
In view of~\eqref{eq2}, we have
\begin{equation}\label{eq3}
R\,(E;\,T)\subseteq
\bigl\{\{E,\,bF,\,c\gamma\}\,\colon\
b\in\mathbb{Z}^+,\ c\in\mathbb{N},\ \gamma\in\{A,C,D\}\bigr\}.
\end{equation}
Note that
\begin{equation}\label{E+F}
E+F=\pi+B.
\end{equation}
So every sum-$2\pi$-collection~$\{E,\,bF,\,c\gamma\}$
with~$b\ge1$ corresponds to a sum-$(\pi-B)$-collection
\begin{equation}\label{eq51}
\{b'F,\,c\gamma\},
\end{equation}
where~$b',c\in\mathbb{N}$ and~$\gamma\in\{A,C,D\}$.
By Lemma~\ref{lem-finite} and the limit~\eqref{lim-angle-F},
there exists~$\delta_4<\delta_3$ such that
for any~$x<\delta_4$, we have
\[
b'F+c\gamma\ne\pi-B,\quad
\forall\,b'\in\mathbb{Z}^+,\
c\in\mathbb{N},\
\gamma\in\{A,C,D\}.
\]
Let~$x<\delta_4$.
In view of~\eqref{eq51},
we find~$b'=0$ and thus~$b=1$.
By~\eqref{E+F}, we see that~$c\ge1$ in~\eqref{eq3}.
This proves that
\begin{equation}\label{Re}
R\,(E;\,T)\subseteq
\bigl\{\{E,\,F,\,c\gamma\}\,\colon\
c\in\mathbb{Z}^+,\
\gamma\in\{A,C,D\}\bigr\}.
\end{equation}

In the same vein, we have
\begin{equation}\label{form-F1}
R\,(F;\,T)\subseteq
\bigl\{\{aF,\,b\beta,\,c\gamma\}\,\colon\
a\in\mathbb{Z}^+,\
b,c\in\mathbb{N},\
\beta,\gamma\in\{E,A,C,D\}\bigr\}.
\end{equation}
By Lemma~\ref{lem-finite} and the limit~\eqref{lim-angle-F},
there exists~$\delta_5<\delta_4$ such that
for any~$x<\delta_5$,
\[
aF+b\beta+c\gamma\ne2\pi,\quad
\forall\,a\in\mathbb{Z}^+,\
b,c\in\mathbb{N},\
\beta,\gamma\in\{A,C,D\}.
\]
In view of~\eqref{form-F1}, we deduce that
\[
R\,(F;\,T)\subseteq
\bigl\{\{aF,\,bE,\,c\gamma\}\,\colon\
a,b\in\mathbb{Z}^+,\
c\in\mathbb{N},\
\gamma\in\{A,C,D\}\bigr\}.
\]
So any sum-$2\pi$-collection to~$F$ is a sum-$2\pi$-collection to~$E$.
By~\eqref{Re}, we obtain that
\[
R\,(F;\,T)=R\,(E;\,T).
\]
This proves the condition~\rn{2}.

Let~$\delta_E=\delta_5$,
and let~$E\in AB$ with~$AE<\delta_E$.
By Lemma~\ref{lem-finite},
there exists~$\delta_F$ depending on the choice of~$E$
such that the condition~\rn{3} holds
for any~$F\in BC$ with~$BF<\delta_F$.
This completes the proof.
\qed

We remark that the points~$E$ and~$F$ can be chosen
from any other pair of adjacent edges of~$Q$,
subject to analogous conditions.
This idea will be employed in the proof of Theorem~\ref{thm-Q}.

We say that a cross-section is {\em proper} if
none of its vertices is a vertex of~$\p$.
Let~$\mathcal{C}_\p$ be the set of
proper pentagonal cross-sections of~$\p$.
Fix two points
\begin{align*}
&E^0\in AB,\\[5pt]
&F^0\in BC,
\end{align*}
satisfying the conditions~\rn{1}---\rn{3},
and write the pentagon
\begin{equation}\label{P0}
P^0=A^0E^0F^0C^0D^0,
\end{equation}
where~$A^0=A$, $C^0=C$, and~$D^0=D$.
Now we recursively define a sequence
$\{P^k\}_{\,k\ge1}$
of proper pentagonal cross-sections.

Let $\pi$ be a plane which crosses~$\p$ nontrivially.
Let~$l$~be a line in~$\pi$.
For any~$\e>0$,
denote by~$\pi_+^\e$ (resp.~$\pi_-^\e$) the plane obtained by
rotating~$\pi$ around~$l$ by the angle~$\e$ (resp.~$-\e$).
It is clear that there exists $\e$ such that
at least one of the planes~$\pi_+^\e$ and~$\pi_-^\e$
crosses~$\p$ nontrivially.
Write
\[
p(\pi;\,l;\,\e)=
\begin{cases}
\pi_+^\e,&\textrm{if~$\pi_+^\e$ crosses~$\p$ nontrivially};\\[5pt]
\pi_-^\e,&\textrm{otherwise}.
\end{cases}
\]
Intuitively,
the plane~$p(\pi;\,l;\,\e)$ is obtained
by rotating~$\pi$ a little along~$l$.
For simplifying notation, we use
$\C(\pi;\,l;\,\e)$
to denote the intersection $C(p(\pi;\,l;\,\e))$
of the plane~$p(\pi;\,l;\,\e)$ and the polyhedron~$\p$.

Since every vertex of~$Q$ has valence~$3$,
there exists~$\delta_0>0$ such that
\[
\C(P^0;\,E^0F^0;\,\e)\in\mathcal{C}_\p,\quad
\forall\,0<\e\le\delta_0.
\]
Define
\begin{equation}\label{P1}
P^1=\C(P^0;\,E^0F^0;\,\delta_0)=A^1E^1F^1C^1D^1
\end{equation}
to be the initial cross-section.
In particular,
we have
\begin{align*}
&E^1=E^0,\\[5pt]
&F^1=F^0.
\end{align*}
Suppose that~$P^k\in\mathcal{C}_\p$ is well-defined
for some~$k\ge1$.
Let~$e_1^k$, $e_2^k$, $\ldots$, $e_5^k$ be the edges of~$P^k$.
It is clear that there exists~$0<\delta_k<\delta_{k-1}$ such that
for any~$0<\e\le\delta_k$ and any edge~$e_j^k$ of~$P^k$,
we have
\[
\C(P^k;\,e_j^k;\,\e)\in\mathcal{C}_\p.
\]
Choosing an edge $e_j^k$, we can define
\begin{equation}\label{def-P}
P^{k+1}=\C(P^k;\,e_j^k;\,\delta_k).
\end{equation}
Note that the cross-section~$P^{k+1}$
depends on the choices of~$e_j^k$ and~$\delta_k$;
while the value of~$\delta_k$ depends on~$P^k$ but
is independent of the choice of~$e_j^k$.
We call the above procedure of getting~$P^{k+1}$ from~$P^k$
the {\em slight-rotating operation}.

Since all cross-sections $P^k$ are proper,
the slight-rotating operation has a certain sign-preserving property
if we take $\delta_k$ small enough.
Denote by~$\mathrm{sgn}(x)$ the signum function, i.e.,
for any real number~$x$,
\[
\mathrm{sgn}(x)=\left\{
\begin{array}{ll}
1, & \hbox{if~$x>0$;} \\
0, & \hbox{if~$x=0$;} \\
-1, & \hbox{if~$x<0$.}
\end{array}
\right.
\]
Let~$i\ge0$, $N>0$, $1\le j\le 5$ and
$a_1$, $\ldots$, $a_5$, $b_1$, $\ldots$, $b_5\in\mathbb{N}$.
Let~$V_i^1$, $\ldots$, $V_i^5$ be the angles of~$P^i$.
Let
\begin{align*}
x_i&=\mathrm{sgn}\biggl(N-\sum_{l=1}^5a_lV_i^l\biggr),\\[5pt]
y_i&=\mathrm{sgn}\biggl(e_j^i-\sum_{l=1}^5b_le_l^i\biggr).
\end{align*}
If~$x_iy_i\ne0$, then
there exists~$0<\delta\le\delta_i$ such that
for any $1\le j\le 5$,
the cross-section
\[
P^{i+1}=\C(P^i;\,e_j^i;\,\delta)
\]
satisfies
\begin{align*}
&x_{i+1}=x_i,\\[5pt]
&y_{i+1}=y_i.
\end{align*}

It is easy to show the above property
if one regards $P^{i+1}$ as a continuous function of the variable $\delta_i$
in the definitions~\eqref{P1} and~\eqref{def-P}.
With the aid of this property,
we deduce that
each cross-section in the sequence $\{P^k\}_{k\ge1}$
satisfies some conditions analogous to~\rn{1}---\rn{3}
as if $\delta_i$'s are small enough.
For all $k\ge2$, we name the vertices of the pentagon~$P^k$ by
\[
P^k=A^kE^kF^kC^kD^k
\]
in the natural way that $E^k\in AB$ and $F^k\in BC$.

\begin{thm}\label{thm-Pk}
For any~$k\ge1$, there exists~$\delta_{k-1}^*\le\delta_{k-1}$
such that for any edge~$e_{j_1}^1$ of~$P^1$,
any edge~$e_{j_2}^2$ of~$P^2$, $\ldots$,
and any edge~$e_{j_{k-1}}^{k-1}$ of~$P^{k-1}$,
the cross-section~$P^k$ defined by~\eqref{def-P} satisfies
\begin{itemize}
\item[\rnp{1}]
$2A^k E^k<\min\{E^kF^k,\,F^kC^k,\,C^kD^k,\,D^kA^k\}$;
\item[\rnp{2}]
$R\,(E^k;\,P^k)=R\,(F^k;\,P^k)\subseteq
\bigl\{\{E^k,\,F^k,\,c\gamma\}\,\colon\
c\in\mathbb{Z}^+,\
\gamma\in\{A^k,\,C^k,\,D^k\}\bigr\}$;
\item[\rnp{3}]
for any~$a,b\in\mathbb{N}$, $c\in\{0,\,\pm1,\,\pm2\}$
and~$l\in\{E^kF^k,\,2F^kC^k\}$, there holds
\[
l\ne a\!\cdot\!C^kD^k+b\!\cdot\!D^kA^k+c\!\cdot\!A^kE^k;
\]
\item[\rnp{4}]
if there exist~$a\in\mathbb{Z}^+$, $b,c\in\mathbb{N}$, $N\in\{\pi,2\pi\}$,
and three pairwise distinct angles $V_k^{j_1}$, $V_k^{j_2}$, $V_k^{j_3}$ of $P^k$
such that
\[
N=aV_k^{j_1}+bV_k^{j_2}+cV_k^{j_3},
\]
then for any~$0\le h\le k$, the corresponding angles
$V_h^{j_1}$, $V_h^{j_2}$, $V_h^{j_3}$ of $P^h$ satisfy
\[
N=aV_h^{j_1}+bV_h^{j_2}+cV_h^{j_3}.
\]
\end{itemize}
\end{thm}

Let~$\delta_i^*$ ($i\ge0$) be defined as in the above theorem.

\begin{lem}\label{lem-parallel}
Let~$\p\in\mathscr{E}$ and~$k\in\mathbb{Z}^+$.
Suppose that
\begin{equation}\label{eq4}
\{E^k,\,F^k,\,A^k\}\in R\,(E^k;\,P^k).
\end{equation}
Then we have
\begin{equation}\label{eq6}
\{E^{k+1},\,F^{k+1},\,A^{k+1}\}\not\in R\bigl(E^{k+1};\,\C(P^k;\,C^kD^k;\,\delta_k^*)\bigr).
\end{equation}
Similarly, if
\[
\{E^k,\,F^k,\,C^k\}\in R\,(E^k;\,P^k),
\]
then we have
\begin{equation}\label{eq7}
\{E^{k+1},\,F^{k+1},\,C^{k+1}\}\not\in R\bigl(E^{k+1};\,\C(P^k;\,D^kA^k;\,\delta_k^*)\bigr).
\end{equation}
\end{lem}

\pf By~\eqref{eq4}, we have
\begin{equation}\label{parallel-V34//V51-1}
F^kC^k\parallel D^kA^k.
\end{equation}
Write
\[
\C(P^k;\,C^kD^k;\,\delta_k^*)
=A^{k+1}E^{k+1}F^{k+1}C^kD^k.
\]
If the statement~\eqref{eq6} is false, then we have
\begin{equation}\label{parallel-V34//V51-2}
F^{k+1}C^k\parallel D^kA^{k+1}.
\end{equation}
\begin{center}
\begin{tikzpicture}
\begin{scope}
\draw[thick]
(0,0)node[below=0.5mm]{$A$}coordinate(A)
--(5,0)node[below=0.5mm]{$B$}coordinate(B)
--+(40:4)node[right]{$C$}coordinate(C)
--(60:2.5)node[below right]{$D$}coordinate(D)
--cycle;

\draw[thick]
%the edges AA1,CC1,DD1.
(100:2)coordinate(Atmp)
(C)+(140:2)coordinate(Ctmp)
(D)--+(110:1.6)node[left,yshift=-1mm]{$D^k$}coordinate(D1);

\draw
(2.8,0)node[below]{$E^k$}coordinate(E1)
(B)+(40:1.5)node[right,xshift=-2mm,yshift=-3mm]{$F^k$}coordinate(F1);
%Find the point A1.
\coordinate (Y1) at (intersection of E1--F1 and D--A);
\coordinate (A1) at (intersection of D1--Y1 and A--Atmp);
\node[left of=A1,xshift=6mm]{$A^k$};
%Find the point C1.
\coordinate (X1) at (intersection of E1--F1 and D--C);
\coordinate (C1) at (intersection of D1--X1 and C--Ctmp);
\node[right of=C1]{$C^k$};

\draw[dotted]
(A1)--(E1)--(F1)--(C1)--(D1)--cycle;

\draw(.7,0)node[below]{$E^{k+1}$}coordinate(E2);
%Find the point A2.
\coordinate (Y2) at (intersection of E2--X1 and D--A);
\coordinate (A2) at (intersection of D1--Y2 and A--Atmp);
\node[left of=A2,xshift=4mm]{$A^{k+1}$};
%Find the point F2.
\coordinate (F2) at (intersection of E2--X1 and B--C);
\node[right of=F2,xshift=-6mm,yshift=-1mm]{$F^{k+1}$};

\draw[dashed]
(A2)--(E2)--(F2)--(C1)--(D1)--cycle;

\draw[thick]
(A)--(A1)
(C)--(C1);

\foreach \point in {A1,E1,F1,C1,D1}
\fill[black,opacity=1] (\point) circle (1.5pt);

\foreach \point in {A2,E2,F2}
\fill[black,opacity=1] (\point) circle (1.5pt);
\end{scope}

\begin{scope}[yshift=-14mm]
\node[xshift=36.8mm,text centered]
{Figure 4. The parallel relation~$F^kF^{k+1}C^k\parallel D^kA^kA^{k+1}$.};
\end{scope}
\end{tikzpicture}
\end{center}
It is clear that the points~$F^k$ and~$F^{k+1}$ are distinct,
while the points~$A^k$ and~$A^{k+1}$ are also distinct.
Hence by~\eqref{parallel-V34//V51-1} and~\eqref{parallel-V34//V51-2},
we find parallel facets
\[
F^kF^{k+1}C^k\parallel D^kA^kA^{k+1},
\]
see Figure 4. But $\p\in\mathscr{E}$, a contradiction.
The relation~\eqref{eq7} can be proved similarly.
This completes the proof. \qed

\section{The main result}

In this section,
we determine all universal tiler
by confirming that $\mathscr{E}$ is empty.

Let~$k\in\mathbb{Z}^+$.
Suppose that~$\mathcal{T}_k$ is a tessellation of the plane
by copies of~$P^k$.
Denote the copies used in~$\mathcal{T}_k$ by
\[
\{\,P_i^k=A_i^kE_i^kF_i^kC_i^kD_i^k\,\colon\ i\in\Lambda_k\,\},
\]
where~$\Lambda_k$ is a set.
Note that each copy in $\T^k$ are arranged counter-clockwise
either in the order
\begin{equation}\label{order1}
A_i^k,\ E_i^k,\ F_i^k,\ C_i^k,\ D_i^k,
\end{equation}
or in the order
\[
A_j^k,\ D_j^k,\ C_j^k,\ F_j^k,\ E_j^k.
\]
Denote by~$\mathcal{I}_k$ the set of indices~$i\in\Lambda_k$ such that
the vertices of~$P_i^k$ are arranged counter-clockwise in the order~\eqref{order1}.
Without loss of generality,
we can invariably suppose that~$1\in\mathcal{I}_k$.

A quadrilateral is said to be cyclic
if all its vertices lie on the same circle.

\begin{thm}\label{thm-Q}
The facet~$Q$ is either a parallelogram or a cyclic quadrilateral.
\end{thm}

\pf Suppose to the contrary that $Q$ is neither a parallelogram nor cyclic.
Without loss of generality,
we can suppose that~$C$ is a largest angle of~$Q$.
Consider the cross-section
\[
P^1=\C(P^0;\,E^0F^0;\,\delta_1^*).
\]
For convenience, we rewrite the copies of~$P^1$ as
\[
P_i^1=A_iE_iF_iC_iD_i.
\]
Since~$P^1$ is a tiler, we can express
\[
S(E_1)=\bigl[\,E_1,\,\beta_j,\,X\,\bigr],
\]
where
\[
\beta_j\in\{A_j,\,E_j,\,F_j,\,C_j,\,D_j\},
\]
and~$X$ is a sequence of angles.
By the condition~\rnp{2}, we have
\[
\beta_j\ne E_j.
\]

Assume that
\[
\beta_j=A_j.
\]
By~\rnp{1},
the sequence~$X$ contains no angle~$A_i$.
By~\rnp{2},
we deduce that~
\begin{equation}\label{S-E1}
S(E_1)
=\bigl[\,E_1,\,A_j,\,F_k\,\bigr]
\end{equation}
for some $k\in\Lambda_1$, see Figure 5.
It follows that
\[
E^1+A^1+F^1=2\pi.
\]
By the condition~\rnp{4}, we obtain
\[
C+D=\pi.
\]
So~$Q$ is a parallelogram, a contradiction.
\begin{center}
\begin{tikzpicture}
\begin{scope}
\draw
(0,0)node[above=3.5mm, right=1mm]{$A_1$}coordinate(A1)
--(1.5,0)node[below=3.3mm, right=-1mm]{$F_k$}coordinate(Fk)--(4.5,0.5)
(A1)--(70:1)
(Fk)node[below=3.3mm,left=0.5mm]{$A_j$}--+(250:1)
(1.5,0)node[above=3.3mm,left=-2mm]{$E_1$}coordinate(E1);

\foreach \point in {A1,E1}
\fill[black,opacity=1] (\point) circle (1.5pt);
\end{scope}

\begin{scope}[xshift=5cm]
\draw
(.5,0)--(4,0)
(2.5,0)node[above=3.5mm, right=1mm]{$A_1$}
coordinate(A1)--(2.85,0.95)
(4,0)node[above=3.3mm,left=-2mm]{$E_1$}
coordinate(E1)--(7,.5)
(4,0)node[below=3mm,left=-2mm]{$\beta_j$}--(5,-1)
(2.5,0)node[above=4mm, left=-2mm]{$Y$};

\foreach \point in {A1,E1}
\fill[black,opacity=1] (\point) circle (1.5pt);
\end{scope}

\begin{scope}[xshift=50mm,yshift=-16mm]
\node[xshift=10mm,text centered]
{Figure 5. The arrangements~\eqref{S-E1} and~\eqref{S-A1Y}.};
\end{scope}
\end{tikzpicture}
\end{center}

Below we can suppose that
\[
\beta_j\in\{F_j,\,C_j,\,D_j\}.
\]
In this case, the condition~\rnp{1} implies
\begin{equation}\label{S-A1Y}
S(A_1)=[\,A_1,\,Y,\,\pi\,],
\end{equation}
where~$Y$ is a sequence of angles, see Figure 5.
By~\rnp{1} and~\rnp{2},
we deduce that~$Y$ contains no angle~$A_i$.
Thus there exist~$b\in\mathbb{Z}^+$
and~$\beta^1\in\{C^1\!,\,D^1\}$ such that
\[
A^1+b\beta^1=\pi.
\]
By~\rnp{4}, we find that
\[
A+b\beta=\pi
\]
for some~$\beta\in\{C,D\}$.
Since~$Q$ is neither a parallelogram nor cyclic,
we see that~$b\ge2$. But~$C$ is a largest angle,
so~$\beta=D$. Namely,
\begin{equation}\label{beta}
A+bD=\pi.
\end{equation}

Consider~$E'\in AD$ and~$F'\in CD$
subject to certain conditions
corresponding to~\rnp{1}---\rnp{4}.
Since~$C$ is a largest angle,
we derive that
\begin{equation}\label{beta'}
A+b'B=\pi
\end{equation}
for some~$b'\ge2$.
Adding~\eqref{beta} and~\eqref{beta'} yields
\[
2\pi=2A+bD+b'B\ge2(A+B+D)=2(2\pi-C),
\]
namely~$C\ge\pi$, a contradiction. This completes the proof.
\qed

\begin{thm}\label{thm-cyclic}
The facet $Q$ is cyclic.
\end{thm}

\pf Suppose to the contrary that
$Q$ is a non-cyclic facet.
By Theorem~\ref{thm-Q}, it is a parallelogram.
Without loss of generality,
we can suppose that~$C$ is a smallest angle of~$Q$, and that
\[
A+D=\pi.
\]
Since~$Q$ is non-cyclic,
we deduce that the angles~$C$ and~$D$ have distinct sizes.
Therefore, the condition~\rn{2} implies that
\[
R\,(E^0;\,P^0)
\subseteq\bigl\{\{E^0,\,F^0,\,A^0\},\ \{E^0,\,F^0,\,C^0\}\bigr\}.
\]
Consider the cross-sections (see Figure 6)
\begin{align}
P^1&=\C(P^0;\,E^0F^0;\,\delta_2^*)=A^1E^1F^1C^1D^1,\label{P11}\\[5pt]
P^2&=\C(P^1;\,C^1D^1;\,\delta_2^*)=A^2E^2F^2C^2D^2,\label{P21}\\[5pt]
P^3&=\C(P^2;\,D^2A^2;\,\delta_2^*)=A^3E^3F^3C^3D^3.\label{P31}
\end{align}
\begin{center}
\begin{tikzpicture}
\begin{scope}
%\draw[help lines](-2,-2) grid (10,5);
\draw[thick]
(0,0)node[below]{$A$}coordinate(A)
--(5,0)node[below]{$B$}coordinate(B)
--+(40:4)node[right]{$C$}coordinate(C)
--(60:2.5)node[below right]{$D$}coordinate(D)
--cycle;

\draw[thick]
%the edges AA1,CC1,DD1.
(100:2)coordinate(Atmp)
(C)+(140:2)coordinate(Ctmp)
(D)--+(110:1.6)node[left,yshift=-1mm]{$D^1\!\!=\!\!D^2\!\!=\!\!D^3$}coordinate(D1);

%The cross-section P1 is determined by the plane D1E1F1.
%Pick the points E0=E1 and F0=F1.
\draw
(2.8,0)node[below]{$E^0\!\!=\!\!E^1$}coordinate(E1)
(B)+(40:1.5)node[right,xshift=-2mm,yshift=-3mm]{$F^0\!\!=\!\!F^1$}coordinate(F1);
%Find the point A1.
\coordinate (Y1) at (intersection of E1--F1 and D--A);
\coordinate (A1) at (intersection of D1--Y1 and A--Atmp);
\node[left of=A1,xshift=6mm]{$A^1$};
%Find the point C1.
\coordinate (X1) at (intersection of E1--F1 and D--C);
\coordinate (C1) at (intersection of D1--X1 and C--Ctmp);
\node[right of=C1]{$C^1\!\!=\!\!C^2$};

\draw[dotted]
(A1)--(E1)--(F1)--(C1)--(D1)--cycle;

%The cross-section P2 is determined by the plane C1D1E2.
%Pick the point E2.
\draw(.7,0)node[below]{$E^2$}coordinate(E2);
%Find the point A2.
\coordinate (Y2) at (intersection of E2--X1 and D--A);
\coordinate (A2) at (intersection of D1--Y2 and A--Atmp);
\node[left of=A2,xshift=2mm]{$A^2\!\!=\!\!A^3$};
%Find the point F2.
\coordinate (F2) at (intersection of E2--X1 and B--C);
\node[right of=F2,xshift=-6mm,yshift=-1mm]{$F^2$};

\draw[dashed]
(A2)--(E2)--(F2)--(C1)--(D1)--cycle;

%The cross-section P3 is determined by the plane D1A2E3.
%Pick the point E3.
\draw(1.8,0)node[below]{$E^3$}coordinate(E3);
%Find the point F3.
\coordinate (F3) at (intersection of E3--Y2 and B--C);
\node[right of=F3,xshift=-7mm,yshift=-2mm]{$F^3$};
%Find the point C3.
\coordinate (X3) at (intersection of E3--F3 and D--C);
\coordinate (C3) at (intersection of D1--X3 and C--Ctmp);
\node[right of=C3,xshift=-4mm]{$C^3$};

\draw[thin]
(A2)--(E3)--(F3)--(C3)--(D1)--cycle;

\draw[thick]
(A)--(A1)
(C)--(C3);

\foreach \point in {A1,E1,F1,C1,D1}
\fill[black,opacity=1] (\point) circle (1.5pt);

\foreach \point in {A2,E2,F2}
\fill[black,opacity=1] (\point) circle (1.5pt);

\foreach \point in {E3,F3,C3}
\fill[black,opacity=1] (\point) circle (1.5pt);

\end{scope}

\begin{scope}[yshift=-14mm]
\node[xshift=41mm,text centered]
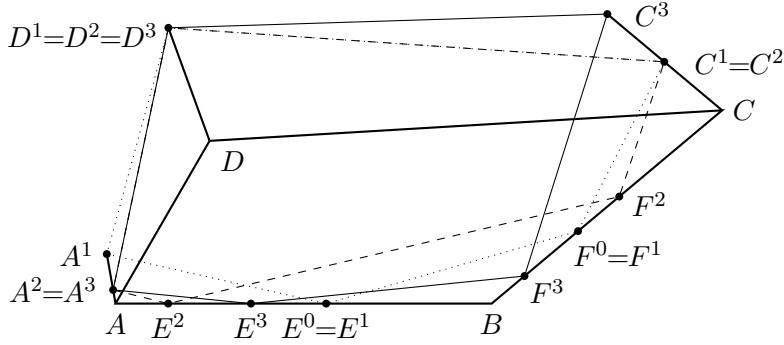
{Figure 6. The cross-sections~$P^1$ (dotted), $P^2$ (dashed) and~$P^3$ (thin).};
\end{scope}
\end{tikzpicture}
\end{center}
By Lemma~\ref{lem-parallel} and~\rnp{4},
we see that
\[
R\,(E^3;\,P^3)=\varnothing.
\]
But the cross-section~$P^3$ tiles the plane, a contradiction.
This completes the proof.
\qed

To proceed further,
we need the following technical lemma.

\begin{lem}\label{lem-irrepresentable}
Let~$k\ge1$. Suppose that
\begin{equation}\label{V124}
R\,(E^k;\,P^k)=R\,(F^k;\,P^k)=\bigl\{\{E^k,\,F^k,\,D^k\}\bigr\}.
\end{equation}
Then any edge~$F_i^kC_i^k$ is not represented.
Moreover, any line segment
\[
F_i^kC_i^k+C_j^kF_j^k
\]
or
\[
E_i^kA_i^k+C_j^kF_j^k
\]
(if exists) is not represented.
\end{lem}

\pf
For convenience,
rewrite the copies~$\{P_i^k\colon i\in\Lambda_k\}$ as
\[
P_i^k=A_iE_iF_iC_iD_i.
\]
Suppose to the contrary that~$F_iC_i$ is represented.
By~\eqref{V124},
there is no point~$F_j$ ($j\ne i$) lying on the
edge~$F_iC_i$, and there is at most one point~$E_j$
lying on~$F_iC_i$.
Therefore, we have
\[
F^kC^k=a_1\!\cdot\!C^kD^k+b_1\!\cdot\!D^kA^k+c_1\!\cdot\!A^kE^k
\]
for some~$a_1,b_1\in\mathbb{N}$ and~$c_1\in\{0,1\}$.
This contradicts to the condition~\rnp{3}.
Hence~$F_iC_i$ is not represented.
For the same reason, any line segment
\[
F_iC_i+C_jF_j
\]
(if exists) is not represented.

Suppose that $1\in\mathcal{I}_k$ and the line segment
\[
E_1F_2=E_1A_1+C_2F_2
\]
is represented in~$\mathcal{T}_k$.
Then~$2\in\mathcal{I}_k$, see Figure 7.
We claim that there exist~$i,j\in\Lambda_k$ such that
\begin{equation}\label{eq1}
S(E_1)=[\,E_1,\,F_i,\,D_j\,].
\end{equation}
\begin{center}
\begin{tikzpicture}
\begin{scope}
\draw
(0,0)node[left=5mm]{$D_j$}--(4.5,0)
(0,0)node[above]{$F_i$}--(160:2)
(3,0)node[above=3.5mm, left=2mm]{$C_i$}coordinate(Ci)--+(120:1)
(3,0)node[above=3.5mm, right=1mm]{$A_h$}--+(70:.9)
(4.5,0)node[above=3.5mm, left=-3mm]{$E_h$}coordinate(Eh)--+(20:1.5)

(0,0)node[below=3mm,right=-2mm]{$E_1$}coordinate(E1)
--(200:2)node[below=3mm, right=-2mm]{$F_1$}coordinate(F1)
--+(230:1)
(1.5,0)node[below=3mm, left=1mm]{$A_1$}
coordinate(A1)--+(250:1)
(1.5,0)node[below=3mm, right=2mm]{$C_2$}--+(-60:1.1)
(4.5,0)node[below=3.5mm, left=-2mm]{$F_2$}--+(-30:1.5);

\foreach \point in {A1,E1,Ci,Eh,F1}
\fill[black,opacity=1] (\point) circle (1.5pt);
\end{scope}

\begin{scope}[yshift=-20mm]
\node[xshift=18mm,text centered]
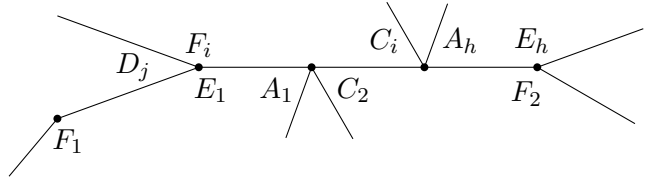
{Figure 7. The tessellation of the~$\e$-neighborhood of the
line segment~$E_1F_2$.};
\end{scope}
\end{tikzpicture}
\end{center}
In fact, by~\eqref{V124},
there is at most one copy~$P_i^k$ such that
the vertex~$E_i$ lying on the line segment~$E_1F_2$,
where~$i\ne1$, $2$.
Also there is at most one point~$F_j$ ($j\ne1,\,2$) lying
on the line segment~$E_1F_2$.
If the sequence~$S(E_1F_2)$ contains an edge~$E_iF_i$,
then the condition~\rnp{2} implies~\eqref{eq1}.
Otherwise, by~\rnp{3},
we can deduce that~$S(E_1F_2)$ contains an edge~$F_iC_i$.
Thus the condition~\rnp{1} yields
\[
S(E_1A_1,\,C_2F_2)=[F_iC_i,\,A_hE_h]
\]
for some~$i$ and~$h$.
In particular, the relation~\eqref{eq1} holds.
The proves the claim.

By~\eqref{eq1},
the edge~$E_1F_1$ is represented.
It follows that
there is at most one point~$E_i$ ($i\ne1$)
lying on the edge~$E_1F_1$,
and there is no point~$F_j$ ($j\ne1$) lying on~$E_1F_1$.
So there exist~$a_2,b_2\in\mathbb{N}$ and~$c_2\in\{0,1\}$ such that
\[
E^kF^k=a_2\!\cdot\!C^kD^k+b_2\!\cdot\!D^kA^k+c_2\!\cdot\!A^kE^k,
\]
contradicting to the condition~\rnp{3}.
This completes the proof. \qed

Here is the main result of this paper.

\begin{thm}\label{main}
A convex polyhedron is a universal tiler if and only if
it is a tetrahedron or a triangular prism.
\end{thm}

\pf It suffices to show that~$\mathscr{E}=\varnothing$.
Suppose to the contrary that~$\p\in\mathscr{E}$.
By Theorem~\ref{thm-cyclic},
we see that~$Q$ is cyclic.
Suppose that~$D$ is a smallest angle of~$Q$.
By the condition~\rn{2}, we have
\begin{equation}\label{eq5}
R\,(E^0;\,P^0)\subseteq\bigg\{
\{E^0,\,F^0,\,A^0\},\
\{E^0,\,F^0,\,C^0\},\
\{E^0,\,F^0,\,D^0\}\bigr\}.
\end{equation}
Consider the cross-sections $P^1$, $P^2$, $P^3$
defined by~\eqref{P11}---\eqref{P31}.
By~\eqref{eq5}, Lemma~\ref{lem-parallel} and
the condition~\rnp{4}, we see that
\[
R\,(E^3;\,P^3)\subseteq\bigl\{\{E^3,\,F^3,\,D^3\}\bigr\}.
\]
Since~$P^3$ is a tiler, by~\rnp{2}, we have
\begin{equation}\label{R-EF3}
R\,(E^3;\,P^3)=R\,(F^3;\,P^3)=\bigl\{\{E^3,\,F^3,\,D^3\}\bigr\}.
\end{equation}
By~\eqref{R-EF3} and~\rnp{4}, we deduce that
\begin{equation}\label{R-Ae}
R_\pi(A^3;\,P^3)\subseteq\bigl\{\{A^3,\,C^3\}\bigr\}\cup
\bigl\{\{bA^3,\,aD^3\}\,\colon\ b\in\mathbb{Z}^+,\,a\in\mathbb{N}\,\bigr\}.
\end{equation}

Consider the tessellations of the plane by copies of~$P^3$.
For convenience, rewrite
\begin{align*}
\mathcal{T}&=\mathcal{T}_3,\\[5pt]
P_i^3&=A_iE_iF_iC_iD_i,\\[5pt]
\Lambda_3&=\Lambda,\\[5pt]
\mathcal{I}_3&=\mathcal{I}.
\end{align*}
By~\eqref{R-EF3} and~\rnp{1},
there exist~$m\in\mathbb{Z}^+$ and~$i_1,i_2,\ldots,i_m\in\Lambda_3$
such that
\begin{equation}\label{form-A1}
S(A_1)=
[\,A_1,\,\alpha_{i_1},\,\alpha_{i_2},\,\ldots,\,\alpha_{i_m},\,\pi\,],
\end{equation}
where
\[
\alpha_{i_j}\in\{A_{i_j},\,C_{i_j},\,D_{i_j}\}.
\]
Let~$\e>0$.
Considering the tessellation of~$\e$-neighborhood of the point~$E_{i_j}$,
we see that
\begin{equation}\label{ineq2}
\alpha_{i_j}\ne A_{i_j},\quad \forall\, 1\le j\le m.
\end{equation}
Write~$i_1=2$.

Assume that~$m=1$.
If~$\alpha_2=D_2$, then
\[
A^3+D^3=\pi
\]
and thus
\[
\{E^3,\,F^3,\,C^3\}\in R\,(E^3;\,P^3),
\]
contradicting to~\eqref{R-EF3}.
Therefore, by~\eqref{ineq2},
the expression~\eqref{form-A1} reduces to
\[
S(A_1)=[\,A_1,\,C_2,\,\pi\,],
\]
see Figure 8.
In this case, if~$2\not\in\mathcal{I}$, then the point~$F_2$ lies
on the line determined by~$D_1A_1$.
By~\eqref{R-EF3},
we find the edge~$F_2C_2$ represented,
contradicting to Lemma~\ref{lem-irrepresentable}.
So~$2\in\mathcal{I}$.
By~\eqref{R-EF3},
the line segment
\[
E_1A_1+C_2F_2
\]
is represented,
also contradicting to Lemma~\ref{lem-irrepresentable}.
\begin{center}
\begin{tikzpicture}
\begin{scope}
\draw
(0,0)node[above=3.5mm,right=-0.8mm]{$A_1$}coordinate(A1)
--++(0.8,0)node[above]{$E_1$}coordinate(E1)
--++(10:1.2)

(0,0)--++(87:2.5)node[below=9mm,right=-2mm]{$D_1$}coordinate(D1)
--+(-62:1.5)
(D1)--++(87:0.5)node[left=0mm]{$F_2$}coordinate(F2)
--+(102:.8)

(-1.5,0)--(A1)node[above=3.5mm,left=-1mm]{$C_2$};

\foreach \point in {A1,E1,D1,F2}
\fill[black,opacity=1] (\point) circle (1.5pt);
\end{scope}

\begin{scope}[xshift=7cm]
\draw
(0,0)node[above=3.5mm,right=-0.8mm]{$A_1$}coordinate(A1)
--++(0.8,0)node[above]{$E_1$}coordinate(E1)
--++(10:1.2)

(0,0)--++(87:2.5)node[below=9mm,right=-2mm]{$D_1$}coordinate(D1)
--+(-62:1.5)

(A1)node[above=3.5mm,left=-1mm]{$C_2$}
--(-3,0)node[above=0.5mm]{$F_2$}coordinate(F2)
--+(160:0.8);

\foreach \point in {A1,E1,D1,F2}
\fill[black,opacity=1] (\point) circle (1.5pt);
\end{scope}

\begin{scope}[yshift=-10mm]
\node[xshift=38mm,text centered]
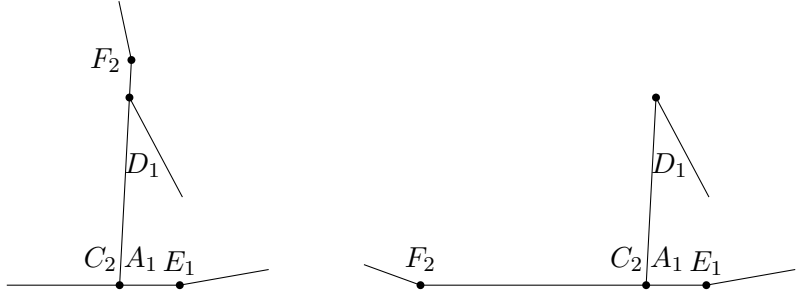
{Figure 8. The tessellation for the case~$m=1$.};
\end{scope}
\end{tikzpicture}
\end{center}

Below we can suppose that~$m\ge2$. Write~$i_2=3$.
By~\eqref{form-A1} and~\eqref{ineq2},
the expression~\eqref{form-A1} reduces to
\begin{equation}\label{S-A1}
S(A_1)=
[\,A_1,\,D_2,\,D_3,\,\ldots,\,\pi\,].
\end{equation}
We have two cases depending on whether~$2\in\mathcal{I}$.

Assume that~$2\in\mathcal{I}$.
By Lemma~\ref{lem-irrepresentable},
the edge~$F_2C_2$ is not represented.
So there exist $s\in\mathbb{Z}^+$
and~$j_1,j_2\ldots,j_s\in\Lambda_3$ such that
\begin{equation}\label{form-C2}
S(C_2)=
[\,C_2,\,\pi,\,\beta_{j_1},\,\beta_{j_2},\,\ldots,\,\beta_{j_s}\,],
\end{equation}
where~$\beta_i$ is an angle of the copy~$T_i$.
Therefore the edge~$D_2C_2$ is represented.
In view of~\eqref{R-EF3}, \eqref{S-A1}, and~\eqref{form-C2},
there is no point~$E_i$ lying on the line segment~$D_2C_2$,
neither is~$F_i$.
Moreover, by~\eqref{R-EF3},
no point~$A_j$ lies in the interior of~$D_2C_2$,
since otherwise the~$\e$-neighborhood of the point~$E_j$
can not be tiled. For the same reason,
no point~$C_j$ lies in the interior of~$D_2C_2$.
Therefore,
\[
S(D_2C_2)\in\{\,[\,D_3A_3\,],\ [\,D_3C_3\,]\,\}.
\]
So~$j_s=3$.
If
\[
S(D_2C_2)=[\,D_3A_3\,],
\]
then either the edge~$A_3E_3$ is represented (when~$s\ge2$),
which is impossible by the condition~\rnp{1} and~\eqref{R-EF3};
or the line segment
\[
E_3A_3+C_2F_2
\]
is represented (when~$s=1$),
which contradicts to Lemma~\ref{lem-irrepresentable}.
See Figure 9.
\begin{center}
\begin{tikzpicture}
\begin{scope}
\draw
(0,0)node[above=3.5mm,right=-0.8mm]{$A_1$}coordinate(A1)
--++(0.8,0)node[above]{$E_1$}coordinate(E1)
--++(10:1.2)

(0,0)node[above=8mm,xshift=-2.8mm]{$D_2$}--+(87:1.4)
(0,0)--(118:2.8)node[below=1.5mm,right=1mm]{$C_2$}coordinate(C2)
(0,0)node[above=4mm,xshift=-6.4mm]{$D_3$}--(149:1.4)

(C2)node[below=4mm,left=-3.2mm]{$A_3$}
--++(-140:0.8)node[below=3mm,right=-2mm]{$E_3$}coordinate(E3)
--+(-130:0.4)
(C2)--++(25:0.9)node[below=3mm,right=-2mm]{$F_2$}coordinate(F2)
--+(10:0.4)
(C2)--+(-155:1.2);

\foreach \point in {A1,E1,C2,F2,E3}
\fill[black,opacity=1] (\point) circle (1.5pt);

\end{scope}
%Case s=1.
\begin{scope}[xshift=6cm]
\draw
(0,0)node[above=3.5mm,right=-0.8mm]{$A_1$}coordinate(A1)
--++(0.8,0)node[above]{$E_1$}coordinate(E1)
--++(10:1.2)

(0,0)node[above=8mm,xshift=-2.8mm]{$D_2$}--+(87:1.4)
(0,0)--(118:2.8)node[below=1.5mm,right=1mm]{$C_2$}coordinate(C2)
(0,0)node[above=4mm,xshift=-6.4mm]{$D_3$}--(149:1.4)

(C2)node[below=4mm,left=-2.5mm]{$A_3$}
(C2)--++(28:0.9)node[below=3mm,right=-1mm]{$F_2$}coordinate(F2)
--+(13:0.4)
(C2)--++(-152:0.8)node[below=3mm,right=-2mm]{$E_3$}coordinate(F3)
--+(-142:0.4);

\foreach \point in {A1,E1,C2,F2,F3}
\fill[black,opacity=1] (\point) circle (1.5pt);
\end{scope}

\begin{scope}[yshift=-10mm]
\node[xshift=28mm,text centered]
{Figure 9. The tessellation for the case~$S(D_2C_2)=[\,D_3A_3\,]$.};
\end{scope}
\end{tikzpicture}
\end{center}
So
\[
S(D_2C_2)=[\,D_3C_3\,].
\]
In this case (see Figure 10),
either~$F_3C_3$ is represented (when~$s\ge2$),
or
\[
F_2C_2+C_3F_3
\]
is represented (when~$s=1$).
Both of them contradicts to Lemma~\ref{lem-irrepresentable}.
\begin{center}
\begin{tikzpicture}
%A=C=87,D=31,E=170,F=165.
%Case s>=2.
\begin{scope}
\draw
(0,0)node[above=3.5mm,right=-0.8mm]{$A_1$}coordinate(A1)
--++(0.8,0)node[above]{$E_1$}coordinate(E1)
--++(10:1.2)

(0,0)node[above=8mm,xshift=-2.8mm]{$D_2$}--+(87:1.4)
(0,0)--(118:2.8)node[below=1.5mm,right=1mm]{$C_2$}coordinate(C2)
(0,0)node[above=4mm,xshift=-6.4mm]{$D_3$}--(149:1.4)

(C2)node[below=5mm,left=-3mm]{$C_3$}
--++(-140:0.9)node[below=3mm,right=-2mm]{$F_3$}coordinate(F3)
--+(-125:0.4)
(C2)--++(25:0.9)node[below=3mm,right=-2mm]{$F_2$}coordinate(F2)
--+(10:0.4)
(C2)--+(-155:1.2);

\foreach \point in {A1,E1,C2,F2,F3}
\fill[black,opacity=1] (\point) circle (1.5pt);

\end{scope}
%Case s=1.
\begin{scope}[xshift=6cm]
\draw
(0,0)node[above=3.5mm,right=-0.8mm]{$A_1$}coordinate(A1)
--++(0.8,0)node[above]{$E_1$}coordinate(E1)
--++(10:1.2)

(0,0)node[above=8mm,xshift=-2.8mm]{$D_2$}--+(87:1.4)
(0,0)--(118:2.8)node[below=1.5mm,right=1mm]{$C_2$}coordinate(C2)
(0,0)node[above=4mm,xshift=-6.4mm]{$D_3$}--(149:1.4)

(C2)node[below=4mm,left=-2.5mm]{$C_3$}
(C2)--++(28:0.9)node[below=3mm,right=-1mm]{$F_2$}coordinate(F2)
--+(13:0.4)
(C2)--++(-152:0.9)node[below=3mm,right=-2mm]{$F_3$}coordinate(F3)
--+(-137:0.4);

\foreach \point in {A1,E1,C2,F2,F3}
\fill[black,opacity=1] (\point) circle (1.5pt);
\end{scope}

\begin{scope}[yshift=-10mm]
\node[xshift=28mm,text centered]
{Figure 10. The tessellation for the case~$S(D_2C_2)=[\,D_3C_3\,]$.};
\end{scope}
\end{tikzpicture}
\end{center}

Hence we have~$2\not\in\mathcal{I}$.
Considering the tessellation of the~$\e$-neighborhood of the point~$E_2$,
we find that
\begin{equation}\label{form-A2}
S(A_2)=
[\,A_2,\,\pi,\,\gamma_{h_1},\,\gamma_{h_2},\,\ldots,\,\gamma_{h_t}\,],
\end{equation}
where~$t\ge1$ and~$\gamma_{h_i}$ is an angle of the copy~$T_{h_i}$.
Therefore~$D_2A_2$ is represented.
For the same reason as in the case~$2\in\mathcal{I}$,
we deduce that
\[
S(D_2A_2)\in\{\,[\,D_3A_3\,],\ [\,D_3C_3\,]\,\}.
\]
If
\[
S(D_2A_2)=[\,D_3A_3\,],
\]
then either~$A_3E_3$ is represented (when~$t\ge2$)
or
\[
E_2A_2+A_3E_3
\]
is represented (when~$t=1$).
See Figure 11.
Both of them are absurd by the condition~\rnp{1}.
\begin{center}
\begin{tikzpicture}
\begin{scope}
\draw
(0,0)node[above=3.5mm,right=-0.8mm]{$A_1$}coordinate(A1)
--++(0.8,0)node[above]{$E_1$}coordinate(E1)
--++(10:1.2)

(0,0)node[above=8mm,xshift=-2.8mm]{$D_2$}--+(87:1.4)
(0,0)--(118:2.8)node[below=1.5mm,right=1mm]{$A_2$}coordinate(A2)
(0,0)node[above=4mm,xshift=-6.4mm]{$D_3$}--(149:1.4)

(A2)node[below=4mm,left=-3.2mm]{$A_3$}
--++(-140:0.8)node[below=3mm,right=-2mm]{$E_3$}coordinate(E3)
--+(-130:0.4)
(A2)--++(25:0.8)node[below=3mm,right=-2mm]{$E_2$}coordinate(E2)
--+(15:0.4)
(A2)--+(-155:1.2);

\foreach \point in {A1,E1,A2,E2,E3}
\fill[black,opacity=1] (\point) circle (1.5pt);

\end{scope}
%Case s=1.
\begin{scope}[xshift=6cm]
\draw
(0,0)node[above=3.5mm,right=-0.8mm]{$A_1$}coordinate(A1)
--++(0.8,0)node[above]{$E_1$}coordinate(E1)
--++(10:1.2)

(0,0)node[above=8mm,xshift=-2.8mm]{$D_2$}--+(87:1.4)
(0,0)--(118:2.8)node[below=1.5mm,right=1mm]{$A_2$}coordinate(A2)
(0,0)node[above=4mm,xshift=-6.4mm]{$D_3$}--(149:1.4)

(A2)node[below=4mm,left=-2.5mm]{$A_3$}
(A2)--++(28:0.8)node[below=3mm,right=-1mm]{$E_2$}coordinate(E2)
--+(18:0.4)
(A2)--++(-152:0.8)node[below=3mm,right=-2mm]{$E_3$}coordinate(F3)
--+(-142:0.4);

\foreach \point in {A1,E1,A2,E2,F3}
\fill[black,opacity=1] (\point) circle (1.5pt);
\end{scope}

\begin{scope}[yshift=-10mm]
\node[xshift=28mm,text centered]
{Figure 11. The tessellation for the case~$S(D_2A_2)=[\,D_3A_3\,]$.};
\end{scope}
\end{tikzpicture}
\end{center}
So
\[
S(D_2A_2)=[\,D_3C_3\,].
\]
In this case, either~$C_3F_3$ is represented (when~$t\ge2$)
or
\[
F_3C_3+A_2E_2
\]
is represented. See Figure 12.
Both of them contradicts to Lemma~\ref{lem-irrepresentable}.
\begin{center}
\begin{tikzpicture}
\begin{scope}
\draw
(0,0)node[above=3.5mm,right=-0.8mm]{$A_1$}coordinate(A1)
--++(0.8,0)node[above]{$E_1$}coordinate(E1)
--++(10:1.2)

(0,0)node[above=8mm,xshift=-2.8mm]{$D_2$}--+(87:1.4)
(0,0)--(118:2.8)node[below=1.5mm,right=1mm]{$A_2$}coordinate(A2)
(0,0)node[above=4mm,xshift=-6.4mm]{$D_3$}--(149:1.4)

(A2)node[below=4mm,left=-3.2mm]{$C_3$}
--++(-140:0.9)node[below=3mm,right=-2mm]{$F_3$}coordinate(F3)
--+(-125:0.4)
(A2)--++(25:0.8)node[below=3mm,right=-2mm]{$E_2$}coordinate(E2)
--+(15:0.4)
(A2)--+(-155:1.2);

\foreach \point in {A1,E1,A2,E2,F3}
\fill[black,opacity=1] (\point) circle (1.5pt);

\end{scope}
%Case s=1.
\begin{scope}[xshift=6cm]
\draw
(0,0)node[above=3.5mm,right=-0.8mm]{$A_1$}coordinate(A1)
--++(0.8,0)node[above]{$E_1$}coordinate(E1)
--++(10:1.2)

(0,0)node[above=8mm,xshift=-2.8mm]{$D_2$}--+(87:1.4)
(0,0)--(118:2.8)node[below=1.5mm,right=1mm]{$A_2$}coordinate(A2)
(0,0)node[above=4mm,xshift=-6.4mm]{$D_3$}--(149:1.4)

(A2)node[below=4mm,left=-2.5mm]{$C_3$}
(A2)--++(28:0.8)node[below=3mm,right=-1mm]{$E_2$}coordinate(E2)
--+(18:0.4)
(A2)--++(-152:0.9)node[below=3mm,right=-2mm]{$F_3$}coordinate(F3)
--+(-137:0.4);

\foreach \point in {A1,E1,A2,E2,F3}
\fill[black,opacity=1] (\point) circle (1.5pt);
\end{scope}

\begin{scope}[yshift=-10mm]
\node[xshift=28mm,text centered]
{Figure 12. The tessellation for the case~$S(D_2A_2)=[\,D_3C_3\,]$.};
\end{scope}
\end{tikzpicture}
\end{center}

To conclude, the cross-section~$P^3$ does not tile the plane.
This implies that
\[
\mathscr{E}=\varnothing,
\]
and completes the proof.
\qed

\end{document}